\documentclass[11pt, twoside]{article}

\usepackage{amsmath,amssymb,amsthm}

\allowdisplaybreaks[4]
\numberwithin{equation}{section}
\swapnumbers

\theoremstyle{plain}
\newtheorem{thm}{Theorem}[section]
\newtheorem{lem}[thm]{Lemma}
\newtheorem{pro}[thm]{Proposition}

\theoremstyle{definition}
\newtheorem{defsrems}[thm]{Definitions and Remarks}
\newtheorem{nots}[thm]{Notations}
\newtheorem{notsrems}[thm]{Notations and Remarks}
\newtheorem{rem}[thm]{Remark}
\newtheorem{rems}[thm]{Remarks}

\newcommand{\A}{\ensuremath{\mathcal A}}
\newcommand{\C}{\ensuremath{\mathbb C}}
\newcommand{\card}{\ensuremath{\mathrm{card}}}
\newcommand{\converges}{\ensuremath{\xrightarrow{\text{$*$-distr}}}}
\newcommand{\convergesas}{\ensuremath{\xrightarrow{\text{$*$-distr, a.s.}}}}
\newcommand{\dif}{\ensuremath{\,\mathrm d}}
\newcommand{\E}{\ensuremath{\mathbb E}}
\newcommand{\F}{\ensuremath{\mathbb F}}
\newcommand{\Fix}{\ensuremath{\mathrm{Fix}\;}}
\newcommand{\hl}{\rule{.5in}{.1mm}}
\newcommand{\K}{\ensuremath{\mathrm K}}
\renewcommand{\Im}{\ensuremath{\mathrm{Im}\;}}
\newcommand{\M}{\ensuremath{\mathcal M}}
\newcommand{\Mat}{\ensuremath{\mathrm{Mat}\;}}
\newcommand{\MN}{\ensuremath{^{(M,N)}}}
\newcommand{\mn}{\ensuremath{_{m,n}}}
\newcommand{\MNq}{\ensuremath{^{(M_q,N_q)}}}
\newcommand{\mO}{\ensuremath{\mathrm O}}
\newcommand{\N}{\ensuremath{^{(N)}}}
\newcommand{\NC}{\ensuremath{\mathcal{NC}}}
\newcommand{\Nq}{\ensuremath{^{(N_q)}}}
\renewcommand{\P}{\ensuremath{\mathcal P}}

\renewcommand{\Re}{\ensuremath{\mathrm{Re}\;}}
\renewcommand{\S}{\ensuremath{\mathcal S}}
\newcommand{\tr}{\ensuremath{\mathrm{tr}}}
\newcommand{\Z}{\ensuremath{\mathbb Z}}

\begin{document}
%%%
\title{\textbf{Asymptotic Freeness of Random Permutation Matrices from Gaussian Matrices}}
\author{\textbf{Mihail Neagu}}
\date{}
\maketitle
%%%
\begin{abstract}
\noindent We show that an independent family of uniformly distributed random permutation matrices is asymptotically $*$-free from an independent family of square complex Gaussian matrices and from an independent family of complex Wishart matrices, and that in both cases the convergence in $*$-distribution actually holds almost surely.

\noindent An immediate consequence is that, if the rows of a GUE matrix are randomly permuted, then the resulting (non self-adjoint) random matrix has a $*$-distribution which is asymptotically circular; similarly, a random permutation of the rows of a complex Wishart matrix results in a random matrix which is asymptotically $*$-distributed like an R-diagonal element from free probability theory.
\end{abstract}

%%%
\section{Introduction.}
The relation between free probability and random matrices springs to a good extent from the fact that, in several important situations, families of $N\times N$ matrices with independent entries turn out to be asymptotically free as $N\to\infty$. As shown by Voiculescu in \cite{V}, this happens for square matrices with independent complex Gaussian entries, for GUE matrices, and for random unitary matrices. It is also observed in \cite{V} that, besides being asymptotically free from each other, the random matrices in question are asymptotically free from constant diagonal matrices.

A variation on this theme is the case of \emph{complex Wishart matrices}. 
These are $N\times N$ matrices of the form $G^*G$, where $G$ is a random 
$M\times N$ matrix with independent complex Gaussian entries. As shown in 
\cite{CC}, an independent family of complex Wishart matrices is asymptotically free as $M,N\to\infty$ with $M/N\to c$, where $c$ is some fixed positive real number.

Yet another phenomenon of asymptotic freeness which arises in a natural situation is the case of $N\times N$ \emph{uniformly distributed random permutation matrices}. It is shown in \cite{N1} that an independent family of such matrices is asymptotically free as $N\to\infty$.

In this note we investigate extensions of the last case mentioned above. 
In other words, we look for families of random matrices which can be 
adjoined to an independent family of uniformly distributed random 
permutation matrices, so that the resulting family is still asymptotically 
free. In contrast to the situations discussed in \cite{V}, adjoining the 
constant diagonal matrices does not work, although in this case it can be shown that a limit joint
distribution still exists (to see that asymptotic freeness 
cannot occur, it suffices to observe that  every permutation matrix $U$ 
and every diagonal matrix $D$ satisfy the algebraic relation 
$U^*DUD=DU^*DU$). We consider the cases when we adjoin either an 
independent family of square complex Gaussian matrices, or an independent 
family of complex Wishart matrices, and we show that in both situations 
the 
resulting family is asymptotically $*$-free (Theorem \ref{thm1} below), 
and that the convergence in $*$-distribution actually holds almost surely 
(Theorem \ref{thm2} below). The asymptotic freeness results are proved by 
first finding 
explicit formulas for the expectations of the mixed moments of the 
random $N\times N$ matrices 
involved, and then showing that as $N\to\infty$ these formulas converge to 
the appropriate expressions obtained from Speicher's theory of 
non-crossing cumulants; similarly, the almost sure convergence results 
are obtained by first finding explicit formulas for the expectations of 
the variances of the $N\times N$ matrices involved, and then showing that 
these expressions are $\mO(1/N^2)$ as $N\to\infty$.

As a consequence, it is easily obtained that a GUE matrix with randomly permuted rows has asymptotically a circular $*$--distribution, and that a Wishart matrix with randomly permuted rows is asymptotically R-diagonal.

It is also possible to get asymptotic $*$-freeness and almost sure convergence results for the case when we adjoin an independent family of \emph{rectangular} complex Gaussian matrices. This however is not a direct generalization of Theorems \ref{thm1} and \ref{thm2}; for a precise statement, see Remarks \ref{MN1} and \ref{MN2}.

In Section 2 of the paper we set the framework and notations which are necessary for a precise formulation of the results, and in Sections 3 and 4 we state and prove the $*$-distribution and the almost sure convergence theorems.
%%%
\section{Preliminaries.}
%%%
\begin{nots}
(1) Throughout this paper, $s$ is a fixed positive integer and we denote by $\F_s$ the free group on the $s$ (free) generators $g_1,g_2,\ldots,g_s$. The elements of $\F_s$ are ``free words'' in the symbols $g_1,g_1^{-1},g_2,g_2^{-1},\ldots,g_s,g_s^{-1}$. The identity of $F_s$ is the ``empty word'' and will be denoted by $e$.

(2)  If $N\geq1$, we denote by $[N]$ the 
set
\begin{equation}
[N]:=\{1,2,\ldots,N\}
\end{equation}
and by $\S_N$ the set of all permutations of $[N]$. For every permutation $\sigma\in\S_N$, we denote by $\Fix\sigma$ the number of fixed points of $\sigma$:
\begin{equation}
\Fix\sigma:=\card\{i\in[N]\,|\,\sigma(i)=i\}.
\end{equation}
If $w\in\F_s$ and $\sigma_1,\sigma_2,\ldots,\sigma_s\in\S_N$, we denote by $w(\sigma_1,\sigma_2,\ldots,\sigma_s)$ the permutation obtained by replacing each $g_r$ by $\sigma_r$ in the expression for $w$. (If $w=e$, then $w(\sigma_1,\sigma_2,\ldots,\sigma_s)$ is the identity permutation of $\S_N$.)

(3) If $M,N\geq1$, we denote by $\M_{M\times N}$ the set of all $M\times N$ matrices with complex entries.  For an $N\times N$ matrix $A$, we denote by $\tr\N A$ the normalized trace of $A$:
\begin{equation}
\tr\N A:=\frac{1}{N}\sum_{i=1}^NA_{ii}.
\end{equation}
Given a permutation $\sigma\in\S_N$, we denote by $\Mat\sigma$ the $N\times N$ matrix whose $(i,j)$-entry is
\begin{equation}   
(\Mat\sigma)_{ij}:=\begin{cases}1&\text{if $\sigma(j)=i$}\\0&\text{if $\sigma(j)\neq i$.}\end{cases}
\end{equation}
A \emph{permutation matrix} is a matrix which is of the form $\Mat\sigma$ for some $\sigma\in\S_N$. We denote by $\P_N$ the set of all $N\times N$ permutation matrices.
\end{nots}
%%%
\begin{nots}
We work with a fixed probability space $(\Omega,\mathcal F,P)$, over which we will consider \emph{random variables} and \emph{random $M\times N$ matrices} (i.~e., measurable functions $f:\Omega\to\C$ and $A:\Omega\to\M_{M\times N}$). For an integrable random variable $f$ over $(\Omega,\mathcal F,P)$, we denote by $\E(f)$ the expectation of $f$:
\begin{equation}
\E(f):=\int_\Omega f(\omega)\dif P(\omega).
\end{equation}
\end{nots}
%%%
We will need the following version of a statement known as ``Wick's Lemma'' (see \cite{J}).
%%%
\begin{lem}[Wick]\label{Wick}
Let $(f_\lambda)_{\lambda\in\Lambda}$ be an independent family of \emph{complex standard Gaussian} random variables over $(\Omega,\mathcal F,P)$; i.~e., $(\Re f_\lambda,\Im 
f_\lambda)_{\lambda\in\Lambda}$ is an independent family of real Gaussian 
random variables over $(\Omega,\mathcal F,P)$ with mean 0 and variance 
$1/2$.

(1) If $m,n\geq1$ with $m\neq n$, and $\alpha_1,\alpha_2,\ldots,\alpha_m,\beta_1,\beta_2,\ldots,\beta_n\in\Lambda$, then
\begin{equation}
\E\left(f_{\alpha_1}f_{\alpha_2}\cdots f_{\alpha_m}\overline{f_{\beta_1}f_{\beta_2}\cdots f_{\beta_n}}\right)=0.
\end{equation}

(2) If $n\geq1$ and $\alpha_1,\alpha_2,\ldots,\alpha_n,\beta_1,\beta_2,\ldots,\beta_n\in\Lambda$, then
\begin{equation}
\E\left(f_{\alpha_1}f_{\alpha_2}\cdots f_{\alpha_n}\overline{f_{\beta_1}f_{\beta_2}\cdots f_{\beta_n}}\right)=\card\{\tau\in\S_n\,|\,\alpha_i=\beta_{\tau(i)}\;\forall i\in[n]\}.
\end{equation}
\end{lem}
%%%
\begin{notsrems}\label{matrices}
(1) We fix an independent family 
$(f_{r;i,j})_{\substack{r\in[s]\\i,j\geq1}}$ of complex standard Gaussian random variables over $(\Omega,\mathcal F,P)$, and for every $M,N\geq1$, we define random $M\times N$ matrices $G_1\MN,G_2\MN,\ldots,G_s\MN:\Omega\to\M_{M\times N}$ by
\begin{equation}
G_r\MN(\omega)=[f_{r;i,j}(\omega)]_{\substack{i\in[M]\\j\in[N]}}.
\end{equation}
For every $r\in[s]$, we denote $G_r\N:=\frac{1}{\sqrt N}G_r^{(N,N)}$ and $W_r\MN:=\frac{1}{N}\left(G_r\MN\right)^*G_r\MN$. The random $N\times N$ matrices $W_r\MN$ are known as \emph{complex Wishart matrices}.

(2) For every $N\geq1$, we fix \emph{uniformly distributed} random
$N\times N$ permutation matrices $U_1\N,U_2\N,\ldots,U_s\N:\Omega\to\P_N$;
in other words, for every $r\in[s]$ and for every $\sigma\in\S_N$,
\begin{equation}\label{uniform}
P\left[\left(U_r\N\right)^{-1}(\Mat\sigma)\right]=\frac{1}{N!}.
\end{equation}
We will consider ``words'' in the random permutation matrices
$U_1\N,U_2\N,\ldots,U_s\N$ and their inverses; if $w\in\F_s$, we denote  
by
$U_w\N$ the random $N\times N$ matrix on $(\Omega,\mathcal F,P)$ obtained
by replacing each $g_r$ by $U_r\N$ in the expression for $w$. (If $w=e$,
then $U_w\N$ is the $N\times N$ identity matrix.)

(3) Finally, we assume that the family $\left(U_r\N,G_r\N\right)_{r=1}^s$ 
is independent for every $N\geq1$, and that the family 
$\left(U_r\N,W_r\MN\right)_{r=1}^s$is independent for every $M,N\geq1$. For our purposes, it suffices to note that these requirements together with \eqref{uniform} imply that, if $g:\left(\M_{N\times N}\right)^s\times\left(\P_N\right)^s
\to\C$ is a function such that for every $\sigma_1,\ldots,\sigma_s\in\S_N$, the random variable
$$\omega\mapsto g\left(G_1\N(\omega),\ldots,G_s\N(\omega),\Mat\sigma_1,\ldots,\Mat\sigma_s\right)$$
is integrable over $(\Omega,\mathcal F,P)$, then
\begin{multline}\label{indep}
\int_\Omega g\left(G_1\N(\omega),\ldots,G_s\N(\omega),U_1\N(\omega),\ldots,U_s\N(\omega)\right)\mathrm dP(\omega)=
\\=\frac{1}{(N!)^s}\sum_{\sigma_1,\ldots,\sigma_s\in\S_N}\int_\Omega g\left(G_1\N(\omega),\ldots,G_s\N(\omega),\Mat\sigma_1,\ldots,\Mat\sigma_s\right)\mathrm dP(\omega)
\end{multline}
\Big(and that an analogous statement holds with each $G_r\N$ replaced by 
$W_r\MN$\Big).
\end{notsrems}
%%%
The remaining of this section is devoted to recalling some notations and facts related to non-crossing permutations and asymptotic $*$-freeness of families of random matrices.
%%%
\begin{notsrems}
Let $m,n\geq1$.

(1) We will use cycle notation for permutations. If $\tau\in\S_n$, the notation $C\in\tau$ will mean that $C$ is one of the cycles in the disjoint cycle decomposition of $\tau$; a fixed point of $\tau$ will be considered as a cycle of length 1. We denote
\begin{equation}
\#(\tau):=\card\{C\,|\,C\in\tau\}.
\end{equation}
To avoid ambiguities when writing down a cycle $(a_1,a_2,\ldots,a_p)$ of length greater than 2, we will always assume that $a_1<\min\limits_{2\leq i\leq p}a_i$.

A permutation $\tau\in\S_n$ is said to be a \emph{pairing} if every 
$C\in\tau$ has 
length 2 (if $n$ is odd, then $\S_n$ contains no pairings). We denote by $\S_n^{(2)}$ the set of all pairings in $\S_n$.

(2) It is a fact (see \cite{B}) that for every $\tau\in\S_n$, we have that
\begin{equation}\label{nc}
\#(\tau)+\#(\tau^{-1}\gamma_n)\leq n+1,
\end{equation}
where $\gamma_n$ is the permutation
\begin{equation}
\gamma_n:=(1,2,\ldots,n)\in\S_n.
\end{equation}
If equality holds in \eqref{nc}, the permutation $\tau$ is said to be \emph{non-crossing}; in this case, the permutation $\tau^{-1}\gamma_n$ (which is also non-crossing) is denoted by $\K(\tau)$ and is called the \emph{Kreweras complement} of $\tau$. We denote by $\NC_n$ the set of all non-crossing permutations in $\S_n$, and by $\NC_n^{(2)}$ the set of all non-crossing pairings in $\S_n$.

(3) A permutation $\tau\in\S_n$ is said to be \emph{parity alternating} if 
$a$ and $\tau(a)$ have opposite parities for every $a\in[n]$.

A permutation $\tau\in\S_n$ is said to be \emph{parity preserving} if $a$ and $\tau(a)$ have the same parity for every $a\in[n]$; in this case, we denote by $\tau_{\mathrm{odd}}$ (respectively, $\tau_{\mathrm{even}}$) the restriction of $\tau$ to the odd (respectively, even) numbers.

It is a fact that, if $n$ is even and $\tau\in\NC_n^{(2)}$, then $\tau$ is 
parity alternating, and hence $\K(\tau)$ is parity preserving. (This follows easily from the more usual way of introducing non-crossing permutations, namely via the concept of \emph{non-crossing partitions} -- see, for instance, \cite{B}.)

(4) A permutation $\tau\in\S_{m+n}$ is said to be \emph{$(m,n)$-connected} if there exist $i$ and $j$ with $1\leq i\leq m$ and $m+1\leq j\leq m+n$ such that $\tau(i)=j$ or $\tau(j)=i$. We denote by $\S\mn$ the set of all $(m,n)$-connected permutations in $\S_{m+n}$. It is a fact (see \cite{MN}) that for every $\tau\in\S\mn$, we have that
\begin{equation}\label{nca}
\#(\tau)+\#(\tau^{-1}\gamma\mn)\leq m+n,
\end{equation}
where $\gamma\mn$ is the permutation
\begin{equation}
\gamma\mn:=(1,2,\ldots,m)(m+1,m+2,\ldots,m+n)\in\S_{m+n}.
\end{equation}
\end{notsrems}
%%%
We will need the following observation.
%%%
\begin{lem}\label{J}
Let $n\geq1$, let $w_1,w_2,\ldots,w_n\in\F_s$, and let $\tau\in\S_n$. For every cycle $C:=(a_1,a_2,\ldots,a_p)\in\tau$, let $w_C\in\F_s$ be the concatenation of the words $w_{a_1},w_{a_2},\ldots,w_{a_p}$:
\begin{equation}
w_C:=w_{a_1}w_{a_2}\cdots w_{a_p}.
\end{equation}
Then for every $N\geq1$ and for every $\sigma_1,\sigma_2,\ldots,\sigma_s\in\S_N$, we have that
\begin{multline}
\card\left\{J:[n]\to[N]\,\left|\right.w_a(\sigma_1,\sigma_2,\ldots,\sigma_s)(J\tau(a))=J(a)\;\forall a\in[n]\right\}=
\\=\prod_{\text{$C\in\tau$}}\Fix w_C(\sigma_1,\sigma_2,\ldots,\sigma_s).
\end{multline}
\end{lem}
%%%
\begin{proof}
Simply note that, if $J:[n]\to[N]$ is a function such that $w_a(\sigma_1,\sigma_2,\ldots,\sigma_s)(J\tau(a))=J(a)\;\forall a\in[n]$, and if $C:=(a_1,a_2,\ldots,a_p)\in\tau$, then the values $J(a_2),\ldots,J(a_p)$ are completely determined by $J(a_1)$, and $J(a_1)$ is a fixed point of the permutation $w_C(\sigma_1,\sigma_2,\ldots\sigma_s)$.
\end{proof}
%%%
\begin{defsrems}
Here we review some free probability concepts (see \cite{VDN}) which are needed in the sequel. We fix a tracial non-commutative $*$-probability space $(\A,\varphi)$; i.~e., \A\ is a unital $*$-algebra and $\varphi:\A\to\C$ is a unital tracial positive linear functional. Elements of \A\ are called \emph{non-commutative random variables} (or simply, \emph{random variables}).

(1) A family $(\A_j)_{j\in J}$ of unital subalgebras of \A\ is said to be \emph{free} if $\varphi(A_1A_2\cdots A_n)=0$ whenever $\varphi(A_i)=0$ for every $i\in[n]$ and $A_i\in\A_{j_i}$, where $j_1,j_2,\ldots,j_n\in J$ are such that $j_1\neq j_2\neq\cdots\neq j_n$.

More generally, a family $(\S_j)_{j\in J}$ of subsets of \A\ is said to be \emph{free} (respectively, \emph{$*$-free}) if the family $(\A_j)_{j\in J}$ is free, where each $A_j$ is the unital subalgebra (respectively, $*$-subalgebra) of \A\ generated by $\S_j$.

(2) Each permutation $\tau\in\NC_n$ gives rise to a multilinear function $\varphi_\tau:\A^n\to\C$ as follows: if $\tau$ consists of $r$ disjoint cycles of the form $C_i=(a_{i,1},a_{i,2},\ldots,a_{i,n_i})\;(1\leq i\leq r)$, then
\begin{equation}
\varphi_\tau(A_1,A_2,\ldots,A_n):=\prod_{i=1}^r\varphi(A_{a_{i,1}}A_{a_{i,2}}\cdots A_{a_{i,n_i}}).
\end{equation}
For instance, if $\tau=(136)(2)(45)\in\NC_6$, then
$$\varphi_\tau(A_1,A_2,A_3,A_4,A_5,A_6)=\varphi(A_1A_3A_6)\varphi(A_2)\varphi(A_4A_5).$$
The \emph{non-crossing cumulants} of the 
space 
$(\A,\varphi)$ are the multilinear functions $k_n:\A^n\to\C\;(n\geq1)$ defined recursively by the formulas
\begin{equation}
\varphi(A_1A_2\cdots A_n)=\sum_{\tau\in\NC_n}k_\tau(A_1,A_2,\ldots,A_n),
\end{equation}
where, if $\tau\in\NC_n$ consists of $r$ cycles of the form $C_i=(a_{i,1},a_{i,2},\ldots,a_{i,n_i})\;(1\leq i\leq r)$, then
\begin{equation}
k_\tau(A_1,A_2,\ldots,A_n):=\prod_{i=1}^r k_{n_i}(A_{a_{i,1}},A_{a_{i,2}},\ldots,A_{a_{i,n_i}}).
\end{equation}

(3) A random variable $U\in\A$ is said to be a \emph{Haar unitary} if 
$U^*U=UU^*=I$ and $\varphi(U^n)=0$ for every $n\in\Z\setminus\{0\}$.

A random variable $G\in\A$ is  said to be \emph{standard circular} if the only non-zero cumulants involving $G$ and $G^*$ are $k_2(G,G^*)=k_2(G^*,G)=1$.

If $c$ is a positive real number, a random variable $W\in\A$ is said to be \emph{free Poisson of parameter $c$} if $W^*=W$ and $k_n(W,W,\ldots,W)=c$ for every $n\geq1$.

There are deep connections between the non-crossing cumulants and the 
concept of freeness, which underlie the combinatorial description of free 
probability theory (see \cite{S}).  For our purposes, it suffices to 
record 
the following two facts, which can be found, for instance, in \cite{NS1}.

(4) If $(A_j)_{j\in J}$ is a $*$-free family of random variables in 
$(\A,\varphi)$, 
if $\varepsilon_1,\varepsilon_2,\ldots,\varepsilon_n\in\{1,*\}$, and if 
$j_1,j_2,\ldots,j_n\in J$ are such that $\card\{j_i\,|\,i\in[n]\}\geq2$, then
\begin{equation}\label{mixed}
k_n(A_{j_1}^{\varepsilon_1},A_{j_2}^{\varepsilon_2},\ldots,A_{j_n}^{\varepsilon_n})=0.
\end{equation}

(5) If $\{A_1,A_2,\ldots,A_n\}$ and $\{B_1,B_2,\ldots,B_n\}$ are subsets of \A\ which are $*$-free from each other, then
\begin{equation}\label{mom-cum}
\varphi(A_1B_1A_2B_2\cdots A_nB_n)=\sum_{\tau\in\NC_n}k_\tau(A_1,A_2,\ldots,A_n)\varphi_{\K(\tau)}(B_1,B_2,\ldots,B_n),
\end{equation}
where for each $\tau\in\NC_n$, $\K(\tau)$ is the Kreweras complement of $\tau$.
\end{defsrems}
%%%
\begin{notsrems}
We fix a positive real number $c$ and a family $(U_r,G_r,W_r)_{r=1}^s$ of 
random variables in $(\A,\varphi)$ such that each $U_r$ is a Haar unitary, 
each $G_r$ is standard circular, each $W_r$ is free Poisson of parameter 
$c$, and the families $(U_r,G_r)_{r=1}^s$ and 
$(U_r,W_r)_{r=1}^s$ 
are $*$-free.

We will consider ``words'' in the random variables $U_1,U_2,\ldots,U_s$ and their inverses; if $w\in\F_s$, we denote by $U_w$ the random variable obtained by replacing each $g_r$  by $U_r$ in the expression for $w$. (If $w=e$, then $U_w$ is the identity of the algebra \A.) It is immediate from the definition of freeness that for every $w\in\F_s$, we have that
\begin{equation}\label{U}
\varphi(U_w)=\begin{cases}1&\text{if $w=e$}\\0&\text{if $w\neq e$.}\end{cases}
\end{equation}
Also, for every $n\geq1$, for every $r_1,r_2,\ldots,r_n\in[s]$, for every $\varepsilon_1,\varepsilon_2,\ldots,\varepsilon_n\in\{1,*\}$, and for every $w_1,w_2,\ldots,w_n\in\F_s$, we have that
\begin{equation}\label{GU}
\varphi\left(G_{r_1}^{\varepsilon_1}U_{w_1}G_{r_2}^{\varepsilon_2}U_{w_2}\cdots G_{r_n}^{\varepsilon_n}U_{w_n}\right)
=\card\left\{\tau\in\NC_n^{(2)}\,\left|\,\substack{r_a=r_{\tau(a)}\text{ and }\varepsilon_a\neq\varepsilon_{\tau(a)}\;\forall a\in[n],\text{ and}
                                                                                  \\w_{a_1}w_{a_2}\cdots w_{a_p}=e\;\forall(a_1,a_2,\ldots,a_p)\in\K(\tau)}\right\}\right.
\end{equation}
and
\begin{equation}\label{WU}
\varphi\left(W_{r_1}U_{w_1}W_{r_2}U_{w_2}\cdots W_{r_n}U_{w_n}\right)
=\hspace{-1cm}\sum_{\substack{\tau\in\NC_n\text{ such that}\\r_a=r_{\tau(a)}\;\forall a\in[n]\text{ and}\\w_{a_1}w_{a_2}\cdots w_{a_p}=e\;\forall(a_1,a_2,\ldots,a_p)\in\K(\tau)}}
   \hspace{-1cm}c^{\#(\tau)}.
\end{equation}
Equations \eqref{GU}  and \eqref{WU} are obtained by using \eqref{mom-cum}, together with \eqref{U} and the explicit description for the cumulants of $G_1,G_1^*,G_2,G_2^*,\ldots,G_s,G_s^*$ (respectively, $W_1,W_2,\ldots,W_s$).
\end{notsrems}
%%%
\begin{defsrems}
If $\alpha$ is a monomial in the random variables $U_r,G_r,W_r$ and their adjoints, then for every $M,N\geq1$ we denote by $\alpha\MN$ the $N\times N$ random matrix over $(\Omega,\mathcal F,P)$ obtained by replacing each $U_r$ by $U_r\N$, each $G_r$ by $G_r\N$, and each $W_r$ by $W_r\MN$ in the expression for $\alpha$, where $U_r\N,G_r\N$, and $W_r\MN$ are as in Notations \ref{matrices}; if $\alpha$ doesn't involve any $W_r$, we simply write $\alpha\N$ instead of $\alpha\MN$.

Let $\mathcal F$ be a subfamily of of $(U_r,G_r,W_r)_{r=1}^s$, and for every $M,N\geq1$, let $\mathcal F\MN$ be the corresponding subfamily of $\left(U_r\N,G_r\N,W_r\MN\right)_{r=1}^s$.

We say that the subfamilies $\mathcal F\MN$ \emph{converge in 
$*$-distribution to $\mathcal F$} as $M,N\to\infty$ with $M/N\to 
c$, written $\mathcal F\MN\converges\mathcal F$, if
\begin{equation}
\lim_{q\to\infty}\E\left(\tr\Nq\alpha\MNq\right)=\varphi(\alpha)
\end{equation}
for every monomial $\alpha$ in the elements of $\mathcal F$ and their adjoints and for every pair $(M_q)_{q=1}^\infty$ and $(N_q)_{q=1}^\infty$ of increasing sequences with $\lim\limits_{q\to\infty}\frac{M_q}{N_q}=c$.

If $\mathcal F\MN\converges\mathcal F$ and the subfamily $\mathcal F$ is 
$*$-free, we say that the families $\mathcal F\MN$ are 
\emph{asymptotically $*$-free} (as $M,N\to\infty$ with $M/N\to c$).

We say that the subfamilies $\mathcal F\MN$ \emph{converge in 
$*$-distribution almost 
surely to $\mathcal F$} as $M,N\to\infty$ with $M/N\to c$, written 
$\mathcal 
F\MN\convergesas\mathcal F$, if they satisfy the stronger condition that
\begin{equation}
\lim_{q\to\infty}\tr\Nq\alpha\MNq=\varphi(\alpha)\quad\text{almost surely}
\end{equation}
for every monomial $\alpha$ in the 
elements of 
$\mathcal F$ and their adjoints and for every pair $(M_q)_{q=1}^\infty$ and $(N_q)_{q=1}^\infty$ of increasing sequences with $\lim\limits_{q\to\infty}\frac{M_q}{N_q}=c$.

The following facts are known:
\begin{alignat}{4}
\left(U_r\N\right)_{r=1}^s&\converges&&(U_r)_{r=1}^s.\tag{see \cite{N1}}
\\\left(G_r\N\right)_{r=1}^s&\convergesas&&(G_r)_{r=1}^s.\tag{see \cite{V} and \cite{HP}}
\\\left(W_r\MN\right)_{r=1}^s&\convergesas&&(W_r)_{r=1}^s.\tag{see \cite{CC}}
\end{alignat}
In the next two sections we will show that $\left(U_r\N,G_r\N\right)_{r=1}^s\convergesas(U_r,G_r)_{r=1}^s$ and that $\left(U_r\N,W_r\MN\right)_{r=1}^s\convergesas(U_r,W_r)_{r=1}^s$.
\end{defsrems}
The following lemma is crucial to our computations.
%%%
\begin{lem}\label{Nica}
(1) If $e\neq w\in\F_s$, then
\begin{equation}\label{lim1}
\limsup_{N\to\infty}\frac{1}{(N!)^s}\sum_{\sigma_1,\ldots,\sigma_s\in\S_N}\Fix w(\sigma_1,\ldots,\sigma_s)<\infty
\end{equation}
and
\begin{equation}\label{lim2}
\limsup_{N\to\infty}\frac{1}{(N!)^s}\sum_{\sigma_1,\ldots,\sigma_s\in\S_N}\big(\Fix w(\sigma_1,\ldots,\sigma_s)\big)^2<\infty.
\end{equation}
(2) If $e\neq w_1\in\F_s$ and $e\neq w_2\in\F_s$, then
\begin{equation}\label{w1w2}
\limsup_{N\to\infty}\frac{1}{(N!)^s}\sum_{\sigma_1,\ldots,\sigma_s\in\S_N}\Fix w_1(\sigma_1,\ldots,\sigma_s)\cdot\Fix w_2(\sigma_1,\ldots,\sigma_s)<\infty.
\end{equation}
\end{lem}
%%%
\begin{proof}
The proof of \eqref{lim1} is embedded in the proof of Theorem 4.1 in \cite {N1}.

For the proof of \eqref{lim2}, note that for every $N\geq1$, we have that
\begin{align*}
\sum_{\sigma_1,\ldots,\sigma_s\in\S_N}\big(\Fix w(\sigma_1,\ldots,\sigma_s)\big)^2
&=\sum_{i,j=1}^N\card\left\{(\sigma_1,\ldots,\sigma_s)\in\left(\S_N\right)^s\,\left|\,\substack{w(\sigma_1,\ldots,\sigma_s)(i)=i\text{ and}\\w(\sigma_1,\ldots,\sigma_s)(j)=j}\right\}\right.=
\\&=N\cdot\card\left\{(\sigma_1,\ldots,\sigma_s)\in\left(\S_N\right)^s\,\left|\,w(\sigma_1,\ldots,\sigma_s)(1)=1\right\}\right.+
\\&+N(N-1)\cdot\card\left\{(\sigma_1,\ldots,\sigma_s)\in\left(\S_N\right)^s\,\left|\,\substack{w(\sigma_1,\ldots,\sigma_s)(1)=1\text{ and}\\w(\sigma_1,\ldots,\sigma_s)(2)=2}\right\}\right..
\end{align*}
Now \eqref{lim2} follows by a direct application of Proposition 3.1.1 in \cite{N2}.

Finally, \eqref{w1w2} follows from \eqref{lim2} together with a simple application of the Cauchy-Schwarz inequality.
\end{proof}
%%%
As shown in \cite{N1}, \eqref{lim1} immediately implies that $\left(U_r\N\right)_{r=1}^s\converges\,(U_r)_{r=1}^s$, since for every $w\in\F_s$ and for every $N\geq1$, we have that
\begin{equation}
\E\left(\tr\N U_w\N\right)=\frac{1}{N\cdot(N!)^s}\sum_{\sigma_1,\ldots,\sigma_s\in\S_N}\Fix w(\sigma_1,\ldots,\sigma_s).
\end{equation}
We observe in Proposition \ref{Uas} below that this conclusion can be 
strengthened, via the use of \eqref{lim2}, to 
$\left(U_r\N\right)_{r=1}^s\convergesas(U_r)_{r=1}^s$.
%%%
\section{Asymptotic $*$-Freeness Results.}
We work within the framework of Section 2, with the notations established there. The main result is the following:
%%%
\begin{thm}\label{thm1}
$\left(U_r\N,G_r\N\right)_{r=1}^s\hspace{-7pt}\converges(U_r,G_r)_{r=1}^s$ and $\left(U_r\N,W_r\MN\right)_{r=1}^s\hspace{-7pt}\converges(U_r,W_r)_{r=1}^s$.

In particular, the families $\left(U_r\N,G_r\N\right)_{r=1}^s$ and $\left(U_r\N,W_r\MN\right)_{r=1}^s$ are asymptotically $*$-free.
\end{thm}
%%%
Towards the proof of Theorem \ref{thm1}, consider an arbitrary monomial in $U_r$ and $G_r$ (respectively, $U_r$ and $W_r$). Up to a cyclic permutation of the symbols (which does not affect the value of the trace $\varphi$), such a monomial must have the form
\begin{equation}\label{ab}
\alpha:=G_{r_1}^{\varepsilon_1}U_{w_1}G_{r_2}^{\varepsilon_2}U_{w_2}\cdots G_{r_n}^{\varepsilon_n}U_{w_n}
\qquad(\text{respectively, }
\beta:=W_{r_1}U_{w_1}W_{r_2}U_{w_2}\cdots W_{r_n}U_{w_n}),
\end{equation}
where $n\geq1$, $r_1,r_2,\ldots,r_n\in[s]$, $\varepsilon_1,\varepsilon_2,\ldots,\varepsilon_n\in\{1,*\}$, and $w_1,w_2,\ldots,w_n\in\F_s$.

For every $\tau\in\S_n$ and for every cycle $C:=(a_1,a_2,\ldots,a_p)\in\tau$, let $w_C\in\F_s$ be the concatenation of the words $w_{a_1},w_{a_2},\ldots,w_{a_p}$:
\begin{equation}
w_C:=w_{a_1}w_{a_2}\cdots w_{a_p}.
\end{equation}
In view of \eqref{GU} and \eqref{WU}, Theorem \ref{thm1} is equivalent to 
the statement that for every pair of increasing sequences $(M_q)_{q=1}^\infty$ and $(N_q)_{q=1}^\infty$ with $\lim\limits_{q\to\infty}\frac{M_q}{N_q}=c$, we have that
\begin{equation}\label{eq1}
\lim_{q\to\infty}\E\left(\tr\Nq\alpha\Nq\right)
=\card\left\{\tau\in\NC_n^{(2)}\,\left|\,\substack{r_a=r_{\tau(a)}\text{ and }\varepsilon_a\neq\varepsilon_{\tau(a)}\;\forall a\in[n],\\\text{and }w_C=e\;\forall C\in\K(\tau)}\right\}\right.
\end{equation}
and
\begin{equation}\label{eq2}
\lim_{q\to\infty}\E\left(\tr\Nq\beta\MNq\right)=\sum_{\substack{\tau\in\NC_n\text{ such that}\\r_a=r_{\tau(a)}\;\forall a\in[n]\text{ and}\\w_C=e\;\forall C\in\K(\tau)}}
                                                                               c^{\#(\tau)}.
\end{equation}

The proofs of \eqref{eq1} and \eqref{eq2} rely on the following lemma, which gives explicit formulas for the expectations involved (before taking limits).
%%%
\begin{lem}\label{lemma1}
Let $\alpha$ and $\beta$ be as in \eqref{ab}.

(1) For every $N\geq1$, we have that
\begin{equation}\label{eq1lemma}
\E\left(\tr\N\alpha\N\right)
=\hspace{-1cm}\sum_{\substack{\tau\in\S_n^{(2)}\text{ such that}\\r_a=r_{\tau(a)}\text{ and }\varepsilon_a\neq\varepsilon_{\tau(a)}\;\forall a\in[n]}}
  \hspace{-1cm}\frac{N^{\#(\tau)}}{N^{n+1}\cdot(N!)^s}\sum_{\sigma_1,\ldots,\sigma_s\in\S_N}\left(\prod_{C\in\tau^{-1}\gamma_n}\Fix w_C(\sigma_1,\ldots,\sigma_s)\right).
\end{equation}
   
(2) For every $M,N\geq1$, we have that
\begin{equation}\label{eq2lemma}
\E\left(\tr\N\beta\MN\right)
=\sum_{\substack{\tau\in\S_n\text{ such that}\\r_a=r_{\tau(a)}\;\forall a\in[n]}}
   \frac{M^{\#(\tau)}}{N^{n+1}\cdot(N!)^s}\sum_{\sigma_1,\ldots,\sigma_s\in\S_N}\left(\prod_{C\in\tau^{-1}\gamma_n}\Fix w_C(\sigma_1,\ldots,\sigma_s)\right).
\end{equation}
\end{lem}
%%%
\begin{proof}
(1) Using \eqref{indep}, we have that
\begin{align*}
&\E\left(\tr\N\alpha\N\right)
=\int_\Omega\frac{1}{N}\sum_{I,J:[n]\to[N]}\left(G_{r_1}\N(\omega)\right)^{\varepsilon_1}_{I(1),J(1)}\left(U_{w_1}\N(\omega)\right)_{J(1),I(2)}\cdots
                                                                              \\&\hspace{5cm}\cdots\left(G_{r_n}\N(\omega)\right)^{\varepsilon_n}_{I(n),J(n)}\left(U_{w_n}\N(\omega)\right)_{J(n),I(1)}\dif P(\omega)=
\\&=\frac{1}{N}\sum_{I,J:[n]\to[N]}\frac{1}{(N!)^s}\sum_{\sigma_1,\ldots,\sigma_s\in\S_N}
       \int_\Omega\left(G_{r_1}\N(\omega)\right)^{\varepsilon_1}_{I(1),J(1)}\Big(\Mat w_1(\sigma_1,\ldots,\sigma_s)\Big)_{J(1),I(2)}\cdots
                             \\&\hspace{5cm}\cdots\left(G_{r_n}\N(\omega)\right)^{\varepsilon_n}_{I(n),J(n)}\Big(\Mat w_n(\sigma_1,\ldots,\sigma_s)\Big)_{J(n),I(1)}\dif P(\omega)=
\\&=\frac{1}{N\cdot(N!)^s}
       \hspace{-1.3cm}\sum_{\substack{I,J:[n]\to[N]\text{ and}\\\sigma_1,\ldots,\sigma_s\in\S_N\text{ such that}\\w_a(\sigma_1,\ldots,\sigma_s)(I\gamma_n(a))=J(a)\;\forall a\in[n]}}
       \hspace{-1cm}\int_\Omega\left(G_{r_1}\N(\omega)\right)^{\varepsilon_1}_{I(1),J(1)}\cdots\left(G_{r_n}\N(\omega)\right)^{\varepsilon_n}_{I(n),J(n)}\dif P(\omega).
\end{align*}
If it is not the case that $n$ is even and that exactly half of $\varepsilon_1,\ldots,\varepsilon_n$ are 1's (and the other half $*$'s), then part (1) of Lemma \ref{Wick} gives that $\E\left(\tr\N\alpha\N\right)=0,$ as desired (since in this case the sum on the right-hand side of \eqref{eq1lemma} is empty). For the rest of the proof, we assume that $n=2k$ and that $[n]=\{c_1,\ldots,c_k,d_1,\ldots,d_k\}$ with $\varepsilon_{c_i}=1$ and $\varepsilon_{d_i}=*$ for every $i\in[k]$. Then
\begin{align*}
&\E\left(\tr\N\alpha\N\right)=
\\&=\frac{1}{N^{k+1}\cdot(N!)^s}\hspace{-1.4cm}\sum_{\substack{I,J:[n]\to[N]\text{ and}\\\sigma_1,\ldots,\sigma_s\in\S_N\text{ such that}
                                                                                                                 \\w_a(\sigma_1,\ldots,\sigma_s)(I\gamma_n(a))=J(a)\;\forall a\in[n]}}
       \hspace{-1.5cm}\E\left(f_{r_{c_1};I(c_1),J(c_1)}\cdots f_{r_{c_k};I(c_k),J(c_k)}\overline{f_{r_{d_1};J(d_1),I(d_1)}\cdots f_{r_{d_k};J(d_k),I(d_k)}}\right)=
\\&=\frac{1}{N^{k+1}\cdot(N!)^s}\hspace{-1.4cm}\sum_{\substack{I,J:[n]\to[N]\text{ and}\\\sigma_1,\ldots,\sigma_s\in\S_N\text{ such that}
                                                                                                                 \\w_a(\sigma_1,\ldots,\sigma_s)(I\gamma_n(a))=J(a)\;\forall a\in[n]}}
       \hspace{-1.5cm}\card\left\{\tau\in\S_k\,\left|\,\substack{r_{c_i}=r_{d_{\tau(i)}},I(c_i)=J(d_{\tau(i)}),\\\text{and }J(c_i)=I(d_{\tau(i)})\;\forall i\in[k]}\right\}\right.
       \quad\text{by part (2) of Lemma \ref{Wick}.}
\end{align*}
Identifying each $\tau\in\S_k$ with the pairing $\left(c_1,d_{\tau(1)}\right)\left(c_2,d_{\tau(2)}\right)\cdots\left(c_k,d_{\tau(k)}\right)\in\S_n^{(2)}$, we get that
\begin{align*}
&\E\left(\tr\N\alpha\N\right)
=\frac{1}{N^{k+1}\cdot(N!)^s}\hspace{-1cm}\sum_{\substack{I,J:[n]\to[N]\text{ and}\\\sigma_1,\ldots,\sigma_s\in\S_N\text{ such that}
                                                                                                            \\w_a(\sigma_1,\ldots,\sigma_s)(I\gamma_n(a))=J(a)\;\forall a\in[n]}}
       \hspace{-1cm}\card\left\{\tau\in\S_n^{(2)}\,\left|\,\substack{\varepsilon_a\neq\varepsilon_{\tau(a)},r_a=r_{\tau(a)},\\\text {and }J(a)=I\tau(a)\;\forall a\in[n]}\right\}\right.=
\\&=\sum_{\substack{\tau\in\S_n^{(2)}\text{ such that}\\r_a=r_{\tau(a)}\text{ and}\\\varepsilon_a\neq\varepsilon_{\tau(a)}\;\forall a\in[n]}}
       \frac{1}{N^{k+1}\cdot(N!)^s}\sum_{\sigma_1,\ldots,\sigma_s\in\S_N}
       \card\left\{(I,J)\,\left|\,\substack{I,J:[n]\to[N]\text{ with } 
J(a)=I\tau(a)\text{ 
and}\\w_a(\sigma_1,\ldots,\sigma_s)(I\gamma_n(a))=J(a)\;\forall 
a\in[n]}\right\}\right.=
\\&=\sum_{\substack{\tau\in\S_n^{(2)}\text{ such that}\\r_a=r_{\tau(a)}\text{ and}\\\varepsilon_a\neq\varepsilon_{\tau(a)}\;\forall a\in[n]}}
       \frac{N^{\#(\tau)}}{N^{n+1}\cdot(N!)^s}\sum_{\sigma_1,\ldots,\sigma_s\in\S_N}
       \card\left\{J:[n]\to[N]\,\left|\,\substack{w_a(\sigma_1,\ldots,\sigma_s)(J\tau^{-1}\gamma_n(a))=J(a)\\\forall a\in[n]}\right\}\right.=
\\&=\sum_{\substack{\tau\in\S_n^{(2)}\text{ such that}\\r_a=r_{\tau(a)}\text{ and}\\\varepsilon_a\neq\varepsilon_{\tau(a)}\;\forall a\in[n]}}
        \frac{N^{\#(\tau)}}{N^{n+1}\cdot(N!)^s}\sum_{\sigma_1,\ldots,\sigma_s\in\S_N}\left(\prod_{C\in\tau^{-1}\gamma_n}\Fix w_C(\sigma_1,\ldots,\sigma_s)\right)
        \quad\text{by Lemma \ref{J}.}
\end{align*}

(2) Using \eqref{indep}, we have that
\begin{align*}
&\E\left(\tr\N\beta\MN\right)
=\int_\Omega\frac{1}{N}\sum_{I,J:[n]\to[N]}\left(W_{r_1}\MN(\omega)\right)_{I(1),J(1)}\left(U_{w_1}\N(\omega)\right)_{J(1),I(2)}\cdots
       \\&\hspace{5cm}\cdots\left(W_{r_n}\MN(\omega)\right)_{I(n),J(n)}\left(U_{w_n}\N(\omega)\right)_{J(n),I(1)}\dif P(\omega)=
\\&=\frac{1}{N}\sum_{I,J:[n]\to[N]}\frac{1}{(N!)^s}\sum_{\sigma_1,\ldots,\sigma_s\in\S_N}
       \int_\Omega\left(W_{r_1}\MN(\omega)\right)_{I(1),J(1)}\Big(\Mat w_1(\sigma_1,\ldots,\sigma_s)\Big)_{J(1),I(2)}\cdots
       \\&\hspace{5cm}\cdots\left(W_{r_n}\MN(\omega)\right)_{I(n),J(n)}\Big(\Mat w_n(\sigma_1,\ldots,\sigma_s)\Big)_{J(n),I(1)}\dif P(\omega)=
\\&=\frac{1}{N^{n+1}\cdot(N!)^s}\hspace{-1cm}\sum_{\substack{I,J:[n]\to[N],\\K:[n]\to[M],\text{ and}\\\sigma_1,\ldots,\sigma_s\in\S_N\text{ such that}
                                                                                                                 \\w_a(\sigma_1,\ldots,\sigma_s)(I\gamma_n(a))=J(a)\;\forall a\in[n]}}                                                                            
       \hspace{-1cm}\int_\Omega\left(G_{r_1}\MN(\omega)\right)^*_{I(1),K(1)}\left(G_{r_1}\MN(\omega)\right)_{K(1),J(1)}\cdots
       \\&\hspace{5cm}\cdots\left(G_{r_n}\MN(\omega)\right)^*_{I(n),K(n)}\left(G_{r_n}\MN(\omega)\right)_{K(n),J(n)}\dif P(\omega)=
\\&=\frac{1}{N^{n+1}\cdot(N!)^s}\hspace{-1cm}\sum_{\substack{I,J:[n]\to[N],\\K:[n]\to[M],\text{ and}\\\sigma_1,\ldots,\sigma_s\in\S_N\text{ such that}
                                                                                                                 \\w_a(\sigma_1,\ldots,\sigma_s)(I\gamma_n(a))=J(a)\;\forall a\in[n]}}
       \hspace{-1.5cm}\E\left(f_{r_1;K(1,)J(1)}\cdots f_{r_n;K(n),J(n)}\overline{f_{r_1;K(1),I(1)}\cdots f_{r_n;K(n),I(n)}}\right).
\end{align*}
Using part (2) of Lemma \ref{Wick}, we get that
\begin{align*}
&\E\left(\tr\N\beta\MN\right)
=\frac{1}{N^{n+1}\cdot(N!)^s}\hspace{-1cm}\sum_{\substack{I,J:[n]\to[N],\\K:[n]\to[M],\text{ and}\\\sigma_1,\ldots,\sigma_s\in\S_N\text{ such that}
                                                                                                                 \\w_a(\sigma_1,\ldots,\sigma_s)(I\gamma_n(a))=J(a)\;\forall a\in[n]}}
       \hspace{-1.5cm}\card\left\{\tau\in\S_n\,\left|\,\substack{r_a=r_{\tau(a)},K(a)=K\tau(a),\\\text{and }J(a)=I\tau(a)\;\forall a\in[n]}\right\}\right.=
\\&=\sum_{\substack{\tau\in\S_n\text{ such that}\\r_a=r_{\tau(a)}\;\forall a\in[n]}}\frac{1}{N^{n+1}\cdot(N!)^s}\sum_{\sigma_1,\ldots,\sigma_s\in\S_N}
       \card\left\{(I,J,K)\,\left|\,\substack{I,J:[n]\to[N]\text{ and 
}K:[n]\to[M]\text{ with}\\K(a)=K\tau(a),J(a)=I\tau(a),\text{ and}
                                                                \\w_a(\sigma_1,\ldots,\sigma_s)(I\gamma_n(a))=J(a)\;\forall a\in[n]}\right\}\right..
\end{align*}
The requirement that $K(a)=K\tau(a)$ for every $a\in[n]$ is equivalent to the statement that $K$ is constant on each cycle of $\tau$, so it follows that
\begin{align*}
&\E\left(\tr\N\beta\MN\right)=
\\&=\sum_{\substack{\tau\in\S_n\text{ such that}\\r_a=r_{\tau(a)}\;\forall a\in[n]}}\frac{M^{\#(\tau)}}{N^{n+1}\cdot(N!)^s}\sum_{\sigma_1,\ldots,\sigma_s\in\S_N}
    \card\left\{J:[n]\to[N]\,\left|\,\substack{\\w_a(\sigma_1,\ldots,\sigma_s)(J\tau^{-1}\gamma_n(a))=J(a)\\\forall a\in[n]}\right\}\right.=
\\&=\sum_{\substack{\tau\in\S_n\text{ such that}\\r_a=r_{\tau(a)}\;\forall a\in[n]}}
       \frac{M^{\#(\tau)}}{N^{n+1}\cdot(N!)^s}\sum_{\sigma_1,\ldots,\sigma_s\in\S_N}\left(\prod_{C\in\tau^{-1}\gamma_n}\Fix w_C(\sigma_1,\ldots,\sigma_s)\right)
       \quad\text{by Lemma \ref{J}.}
\end{align*}
\end{proof}
%%%
\begin{proof}[Proof of Theorem \ref{thm1}.]
We  only give the proof of \eqref{eq2}, since that of \eqref{eq1} is virtually identical.

Let $(M_q)_{q=1}^\infty$ and $(N_q)_{q=1}^\infty$ be two increasing 
sequences 
with $\lim\limits_{q\to\infty}\frac{M_q}{N_q}=c$. If $\tau\in\S_n$ is such that $\tau^{-1}\gamma_n$ has a cycle $C_0$ with $w_{C_0}\neq e$, then
\begin{align*}
\frac{M_q^{\#(\tau)}}{N_q^{n+1}\cdot(N_q!)^s}&\sum_{\sigma_1,\ldots,\sigma_s\in\S_{N_q}}
                                                                                    \left(\prod_{\text{cycles $C$ of $\tau^{-1}\gamma_n$}}\Fix w_C(\sigma_1,\ldots,\sigma_s)\right)\leq
\\&\leq\frac{M_q^{\#(\tau)}}{N_q^{n+1}\cdot(N_q!)^s}\sum_{\sigma_1,\ldots,\sigma_s\in\S_{N_q}}N_q^{\#(\tau^{-1}\gamma_n)-1}\cdot\Fix w_{C_0}(\sigma_1,\ldots,\sigma_s)=
\\&=\left(\frac{M_q}{N_q}\right)^{\#(\tau)}\cdot\frac{N_q^{\#(\tau)+\#(\tau^{-1}\gamma_n)-(n+1)}}{N_q}\cdot
       \frac{1}{(N_q!)^s}\sum_{\sigma_1,\ldots,\sigma_s\in\S_{N_q}}\Fix w_{C_0}(\sigma_1,\ldots,\sigma_s)
\\&\xrightarrow[q\to\infty]{}0\qquad\text{by \eqref{nc} and \eqref{lim1}.}
\end{align*}
Hence part (2) of Lemma \ref{lemma1} gives that
\begin{align*}
\lim_{q\to\infty}\E\left(\tr\Nq\beta\MNq\right)
&=\sum_{\substack{\tau\in\S_n\text{ such that}\\r_a=r_{\tau(a)}\;\forall a\in[n]\text{ and}\\w_C=e\;\forall C\in\tau^{-1}\gamma_n}}
     \lim_{q\to\infty}\frac{M_q^{\#(\tau)}}{N_q^{n+1}\cdot(N_q!)^s}\sum_{\sigma_1,\ldots,\sigma_s\in\S_{N_q}}N_q^{\#(\tau^{-1}\gamma_n)}=
\\&=\sum_{\substack{\tau\in\S_n\text{ such that}\\r_a=r_{\tau(a)}\;\forall a\in[n]\text{ and}\\w_C=e\;\forall C\in\tau^{-1}\gamma_n}}
       \lim_{q\to\infty}\left(\frac{M_q}{N_q}\right)^{\#(\tau)}\cdot N_q^{\#(\tau)+\#(\tau^{-1}\gamma_n)-(n+1)}
\\&=\sum_{\substack{\tau\in\NC_n\text{ such that}\\r_a=r_{\tau(a)}\;\forall a\in[n]\text{ and}\\w_C=e\;\forall C\in\K(\tau)}}c^{\#(\tau)}
       \qquad\text{by \eqref{nc} and the definition of $\NC_n$.}
\end{align*}
\end{proof}
%%%
\begin{rems}
(1) A corollary of Theorem \ref{thm1} is that the families $\hspace{-1pt}\left(\hspace{-2pt}U_r\N,\frac{G_r\N+\left(G_r\N\hspace{-2pt}\right)^*}{\sqrt2}\right)_{r=1}^s$ are asymptotically $*$-free, which has the meaning that random permutation matrices are asymptotically $*$-free from random self-adjoint matrices with independent Gaussian entries. \Big(The random matrices $\frac{G_r\N+\left(G_r\N\right)^*}{\sqrt2}$ are known as \emph{GUE matrices}.\Big)

(2) In principle, it would be possible to approach the proof of the 
asymptotic 
$*$-freeness of random permutation 
matrices from either Wishart or GUE matrices by attempting to use condition (C) appearing on page 398 of \cite{CC}. In the same vein, towards proving the almost sure convergence results from the next section, one could attempt to use condition (C') appearing on page 415 of \cite{CC}. However, pursuing this alternative method would require some strengthening of the results of 
\cite{N1} and \cite{N2}, since they do not give directly that condition 
(C) 
(or even 
more so, (C')) holds for random permutation matrices. Moreover, this 
approach would not yield the explicit formulas appearing in Lemma 
\ref{lemma1} (and in Lemma \ref{lemma2} in the next section) for the 
expectations of the  
mixed moments (and of the variances) of the $N\times N$ random matrices 
involved. \end{rems}
%%%
\begin{rems}
(1) If $X\N$ is an $N\times N$ GUE matrix with randomly permuted rows, then $X\N$ converges in $*$-distribution to a standard circular random variable as $N\to\infty$.

This is because $X\N$ can be realized as $U\N \cdot\frac{\left(G\N\right)^*+G\N}{\sqrt{2}}$, which by Theorem \ref{thm1} converges in $*$-distribution to $UY$, where $U$ is a Haar unitary, $Y$ is \emph{standard semicircular} (see \cite{V} for the definition), and the family $(U,Y)$ is $*$-free -- but it is known that such a $UY$ is standard circular (see \cite{V}, or \cite{NSS} where this is covered from the point of view of R-diagonal random variables).

(2) If $Y\MN$ is an $N\times N$ complex Wishart matrix with randomly 
permuted rows, then $Y\MN$ converges in $*$-distribution to an 
\emph{R-diagonal} (see \cite{NSS} for the definition) random variable as 
$M,N\to\infty$ with $M/N\to c$.

This is because $Y\MN$ can be realized as $U\N W\MN$, which by Theorem \ref{thm1} converges in $*$-distribution to $UW$, where $U$ is a Haar unitary, $W$ is free Poisson of parameter $c$, and the family $(U,W)$ is $*$-free -- but it is known that such a $UW$ is R-diagonal (see \cite{NSS}).
\end{rems}
%%%
\begin{rem}\label{MN1}
We conclude this section by discussing an analogue of the first statement 
of Theorem \ref{thm1}, 
obtained by allowing the complex Gaussian matrices involved to be 
\emph{rectangular}, of size $M\times N$ with $M/N\to c$.

Let $G_1\MN,\ldots,G_s\MN,U_1\N,\ldots,U_s\N$ be as in Notations 2.4, let 
$T_1^{(M)},\ldots,T_s^{(M)}$ be uniform random $M\times M$ permutation 
matrices such that the family 
$\left(T_r^{(M)},U_r\N,G_r\MN\right)_{r=1}^s$ 
is independent, and consider the $(M+N)\times(M+N)$ random matrices
$$T_r\MN:=\begin{bmatrix}T_r^{(M)}&0\\0&0\end{bmatrix},
\quad U_r\MN:=\begin{bmatrix}0&0\\0&U_r\N\end{bmatrix},
\quad 
H_r\MN:=\frac{1}{\sqrt{M+N}}\begin{bmatrix}0&G_r\MN\\0&0\end{bmatrix}.$$

Let $P$ and $Q$ be projections in \A\ with $P+Q=I$ and 
$\frac{\varphi(P)}{\varphi(Q)}=c$, and let $T_1,\ldots,T_s,U_1,\ldots,U_s$ 
be partial isometries in \A\ such that for every $r\in[s]$, we have that 
$T_r^*T_r=T_rT_r^*=P$, $U_r^*U_r=U_rU_r^*=Q$, and 
$\varphi(T^n)=\varphi(U^n)=0$ for every $n\in\Z\setminus\{0\}$. (Here, if 
$n<0$, $T^n$ 
and $U^n$ are interpreted as $(T^*)^{-n}$ and $(U^*)^{-n}$, respectively.)

As before, if $w\in\F_s$, we denote by $T_w$ (respectively, $U_w$) the element of \A\ obtained by replacing each $g_r$  by $T_r$ (respectively, by $U_r$) in the expression for $w$. (If $w=e$, then $T_w:=P$ and $U_w:=Q$.)

Let $G_1,G_2,\ldots,G_s$ be standard circular elements of \A\ such that 
the 
family $(\{T_r,U_r\},G_r)_{r=1}^s$ is $*$-free, and for every $r\in[s]$, 
let $H_r:=PG_rQ$.

Then
\begin{equation}\label{HTU}
\left(T_r\MN,U_r\MN,H_r\MN\right)_{r=1}^s\converges(T_r,U_r,H_r)_{r=1}^s
\end{equation}
as $M,N\to\infty$ with $M/N\to c$; we sketch the proof below. 

Consider an arbitrary monomial $\alpha$ in $T_r,U_r,H_r$, and their 
adjoints, and let $\alpha\MN$ be the $(M+N)\times(M+N)$ random matrix obtained by replacing each $H_r$ by $H_r\MN$, each $T_r$ by $T_r\MN$, and each $U_r$ by $U_r\MN$ in the expression for $\alpha$. It is easily seen that $\varphi(\alpha)=0=\tr^{(M+N)}\alpha\MN$ unless, up to a cyclic permutation of the symbols (which does not affect the value of the trace $\varphi$), $\alpha$ has the form
\begin{equation}\label{alpha}
\alpha=H_{r_1}^*T_{w_1}H_{r_2}U_{w_2}H_{r_3}^*T_{w_3}H_{r_4}U_{w_4}\cdots 
H_{r_{2k-1}}^*T_{w_{2k-1}}H_{r_{2k}}U_{w_{2k}},
\end{equation}
where $k\geq1$, $r_1,r_2,\ldots,r_{2k}\in[s]$, and $w_1,w_2,\ldots,w_{2k}\in\F_s$.

If $\alpha$ is of the above form, then a very similar argument to the one given in the proof of part (1) of Lemma \ref{lemma1} shows that for every $M,N\geq1$, we have that
\begin{multline}\label{lemma1MN}
\E\left(\tr^{(M+N)}\alpha\MN\right)=\sum_{\substack{\tau\in\S_{2k}^{(2)}\text{ such that}\\r_a=r_{\tau(a)}\;\forall a\in[2k]}}\frac{(M+N)^{\#(\tau)}}{(M+N)^{2k+1}\cdot(M!)^s(N!)^s}\cdot
\\\cdot\sum_{\substack{\sigma_1,\ldots,\sigma_s\in\S_M\\\rho_1,\ldots,\rho_s\in\S_N}}\left(\prod_{C\in(\tau^{-1}\gamma_{2k})_{\mathrm{odd}}}\Fix w_C(\sigma_1,\ldots,\sigma_s)\right)
                                                                                                                                                         \left(\prod_{C\in(\tau^{-1}\gamma_{2k})_{\mathrm{even}}}\Fix w_C(\rho_1,\ldots,\rho_s)\right).
\end{multline}
Let $(M_q)_{q=1}^\infty$ and $(N_q)_{q=1}^\infty$ be two increasing 
sequences 
with $\lim\limits_{q\to\infty}\frac{M_q}{N_q}=c$. If $\tau\in\S_{2k}^{(2)}$ is such that $\tau^{-1}\gamma_{2k}$ has a cycle $C_0$ with $w_{C_0}\neq e$, then a virtually identical argument to the one given in the proof of Theorem \ref{thm1} shows that the term associated to $\tau$ in \eqref{lemma1MN} vanishes in the limit as $q\to\infty$. Thus we have that
\begin{equation}\label{sum}
\lim_{q\to\infty}\E\left(\tr^{(M_q+N_q)}\alpha\MNq\right)
=\hspace{-0.5cm}\sum_{\substack{\tau\in\S_{2k}^{(2)}\text{ such that}\\r_a=r_{\tau(a)}\;\forall a\in[2k]\text{ and}\\w_C=e\;\forall C\in\tau^{-1}\gamma_{2k}}}
   \hspace{-0.5cm}\lim_{q\to\infty}\frac{(M_q+N_q)^{\#(\tau)}M_q^{\#((\tau^{-1}\gamma_{2k})_{\mathrm{odd}})}N_q^{\#((\tau^{-1}\gamma_{2k})_{\mathrm{even}})}}{(M_q+N_q)^{2k+1}}.
\end{equation}
If $\tau\in\S_{2k}^{(2)}\setminus\NC_{2k}^{(2)}$, then the term associated to $\tau$ in \eqref{sum} is zero by \eqref{nc}. Hence
\begin{align}\label{sum1}
\lim_{q\to\infty}\E\left(\tr^{(M_q+N_q)}\alpha\MNq\right)
&=\sum_{\substack{\tau\in\NC_{2k}^{(2)}\text{ such that}\\r_a=r_{\tau(a)}\;\forall a\in[2k]\text{ and}\\w_C=e\;\forall C\in\K(\tau)}}
     \lim_{q\to\infty}\frac{M_q^{\#(\K(\tau)_{\mathrm{odd}})}N_q^{\#(\K(\tau)_{\mathrm{even}})}}{(M_q+N_q)^{\#(\K(\tau))}}=\notag
\\&=\sum_{\substack{\tau\in\NC_{2k}^{(2)}\text{ such that}\\r_a=r_{\tau(a)}\;\forall a\in[2k]\text{ and}\\w_C=e\;\forall C\in\K(\tau)}}\frac{c^{\#(\K(\tau)_{\mathrm{odd}})}}{(1+c)^{\#(\K(\tau))}}.
\end{align}
On the other hand, by replacing each $H_r$ by $PG_rQ$ in \eqref{alpha}, we 
get that
\begin{equation*}
\varphi(\alpha)=\varphi(G_{r_1}^*T_{w_1}G_{r_2}U_{w_2}G_{r_3}^*T_{w_3}G_{r_4}U_{w_4}\cdots 
G_{r_{2k-1}}^*T_{w_{2k-1}}G_{r_{2k}}U_{w_{2k}}),
\end{equation*}
and then an application of \eqref{mom-cum} yields that
\begin{align}\label{sum2}
\varphi(\alpha)
&=\sum_{\substack{\tau\in\NC_{2k}^{(2)}\text{ such that}\\r_a=r_{\tau(a)}\;\forall a\in[2k]\text{ and}\\w_C=e\;\forall C\in\K(\tau)}}
     \varphi(P)^{\#(\K(\tau)_{\mathrm{odd}})}\varphi(Q)^{\#(\K(\tau)_\mathrm{even})}=\notag
\\&=\sum_{\substack{\tau\in\NC_{2k}^{(2)}\text{ such that}\\r_a=r_{\tau(a)}\;\forall a\in[2k]\text{ and}\\w_C=e\;\forall C\in\K(\tau)}}
       \left(\frac{c}{1+c}\right)^{\#(\K(\tau)_{\mathrm{odd}})} \left(\frac{1}{1+c}\right)^{\#(\K(\tau)_{\mathrm{even}})}.
\end{align}       
Now \eqref{HTU} follows from \eqref{sum1} and \eqref{sum2}.
\end{rem}
%%%
\section{Almost Sure Convergence Results.}
We begin this section with a general statement about almost sure convergence, whose proof can be found, for instance, embedded in the proof of Corollary 3.9 in \cite{T}.
%%%
\begin{lem}\label{Steen}
Let $z\in\C$, let $(g_q)_{q=1}^\infty$ be a sequence of integrable random variables on $(\Omega,\mathcal F,P)$ such that $\lim\limits_{q\to\infty}\E(g_q)=z$, and suppose that
\begin{equation}
\sum_{q=1}^\infty\left[\E\left(|g_q|^2\right)-\left|\E(g_q)\right|^2\right]<\infty.
\end{equation}
Then $\lim\limits_{q\to\infty}g_q=z$ almost surely.
\end{lem}
%%%
For the rest of the section, we work within the framework of Section 2, with the notations established there.
%%%
\begin{pro}\label{Uas}
$\left(U_r\N\right)_{r=1}^s\convergesas(U_r)_{r=1}^s$.
\end{pro}
%%%
\begin{proof}
In view of Lemma \ref{Steen}, it suffices to show that for every $w\in\F_s$, we have that
\begin{equation}\label{sumU}
\sum_{N=1}^\infty\left[\E\left(\left|\tr\N U_w\N\right|^2\right)-\left|\E\left(\tr\N U_w\N\right)\right|^2\right]<\infty.
\end{equation}
If $w=e$, then \eqref{sumU} is clear; if $w\neq e$, then
\begin{align*}
&\E\left(\left|\tr\N U_w\N\right|^2\right)-\left|\E\left(\tr\N U_w\N\right)\right|^2=
\\&=\frac{1}{N^2}\cdot\frac{1}{(N!)^s}\sum_{\sigma_1,\ldots,\sigma_s\in\S_N}\big(\Fix w(\sigma_1,\ldots,\sigma_s)\big)^2
  -\frac{1}{N^2}\left(\frac{1}{(N!)^s}\sum_{\sigma_1,\ldots,\sigma_s\in\S_N}\Fix w(\sigma_1,\ldots,\sigma_s)\right)^2=
\\&=\mO(1/N^2)\qquad\text{by part (1) of Lemma \ref{Nica}},
\end{align*}
from which \eqref{sumU} follows.
\end{proof}
%%%
The main result of this section is the following:
%%%
\begin{thm}\label{thm2}
$\left(U_r\N,G_r\N\right)_{r=1}^s\convergesas(U_r,G_r)_{r=1}^s$ and $\left(U_r\N,W_r\MN\right)_{r=1}^s\convergesas(U_r,W_r)_{r=1}^s$.
\end{thm}
%%%
Towards the proof of Theorem \ref{thm2}, consider an arbitrary monomial 
$\alpha$ in $U_r$ and $G_r$ (or in $U_r$ and $W_r$) and their adjoints. In 
view of Lemma 
\ref{Steen}, it suffices to show that, if $(M_q)_{q=1}^\infty$ and $(N_q)_{q=1}^\infty$ are two increasing sequences with $\lim\limits_{q\to\infty}\frac{M_q}{N_q}=c$, then
\begin{equation}
\sum_{q=1}^\infty\left[\E\left(\left|\tr\Nq\alpha\MNq\right|^2\right)-\left|\E\left(\tr\Nq\alpha\MNq\right)\right|^2\right]<\infty.
\end{equation}
In fact, we prove the slightly more general statement that, if $\alpha$ 
and $\beta$ are monomials in $U_r$ and $G_r$ (or in $U_r$ and $W_r$) and 
their adjoints, and 
if $(M_q)_{q=1}^\infty$ and $(N_q)_{q=1}^\infty$ are two increasing sequences with $\lim\limits_{q\to\infty}\frac{M_q}{N_q}=c$, then
\begin{equation}\label{alphabeta}
\E\left(\tr\Nq\alpha\MNq\cdot\tr\Nq\beta\MNq\right)-\E\left(\tr\Nq\alpha\MNq\right)\cdot\E\left(\tr\Nq\beta\MNq\right)=\mO\left(1/N_q^2\right).
\end{equation}
The proof of \eqref{alphabeta} relies on computing its left-hand side explicitly. Whereas Lemma \ref{lemma1} gives an explicit formula for the second term, the following lemma gives an explicit formula for the first term; we omit its proof, since it is virtually identical to that of Lemma \ref{lemma1}.
%%%
\begin{lem}\label{lemma2}
Let $m,n\geq1$, let $r_1,r_2,\ldots,r_{m+n}\in[s]$, let $\varepsilon_1,\varepsilon_2,\ldots,\varepsilon_{m+n}\in\{1,*\}$, and let $w_1,w_2,\ldots,w_{m+n}\in\F_s$.
\item[ (1)] If $\alpha=G_{r_1}^{\varepsilon_1}U_{w_1}G_{r_2}^{\varepsilon_2}U_{w_2}\cdots G_{r_m}^{\varepsilon_m}U_{w_m}$ and $\beta=G_{r_{m+1}}^{\varepsilon_{m+1}}U_{w_{m+1}}G_{r_{m+2}}^{\varepsilon_{m+2}}U_{w_{m+2}}\cdots G_{r_{m+n}}^{\varepsilon_{m+n}}U_{w_{m+n}}$, then for every $N\geq1$ we have that
\begin{multline}
\E\left(\tr\N\alpha\N\cdot\tr\N\beta\N\right)=
\\=\hspace{-0.2cm}\sum_{\substack{\tau\in\S_{m+n}^{(2)}\text{ such that}\\r_a=r_{\tau(a)}\text{ and }\varepsilon_a\neq\varepsilon_{\tau(a)}\;\forall a\in[m+n]}}
     \hspace{-1cm}\frac{N^{\#(\tau)}}{N^{m+n+2}\cdot(N!)^s}\sum_{\sigma_1,\ldots,\sigma_s\in\S_N}
     \left(\prod_{C\in\tau^{-1}\gamma\mn}\Fix w_C(\sigma_1,\ldots,\sigma_s)\right).
\end{multline}
\item[ (2)] If $\alpha=W_{r_1}U_{w_1}W_{r_2}U_{w_2}\cdots W_{r_m}U_{w_m}$ and $\beta=W_{r_{m+1}}U_{w_{m+1}}W_{r_{m+2}}U_{w_{m+2}}\cdots W_{r_{m+n}}U_{w_{m+n}}$, then for every $M,N\geq1$ we have that
\begin{multline}
\E\left(\tr\N\alpha\MN\cdot\tr\N\beta\MN\right)=
\\=\sum_{\substack{\tau\in\S_{m+n}\text{ such that}\\r_a=r_{\tau(a)}\;\forall a\in[m+n]}}
     \frac{M^{\#(\tau)}}{N^{m+n+2}\cdot(N!)^s}\sum_{\sigma_1,\ldots,\sigma_s\in\S_N}\left(\prod_{C\in\tau^{-1}\gamma\mn}\Fix w_C(\sigma_1,\ldots,\sigma_s)\right).
\end{multline}
\end{lem}
%%%
\begin{proof}[Proof of Theorem \ref{thm2}.]
We only present the proof of \eqref{alphabeta} for the case when $\alpha$ 
and $\beta$ are monomials in $U_r,G_r$, and their adjoints, since the 
proof for the case 
when they are monomials in $U_r,W_r$, and their adjoints  is essentially 
the same. For 
simplicity, we will omit the subscript $q$.

Let $\alpha$ and $\beta$ be as in part (1) of Lemma \ref{lemma2}. Writing 
$\sum\limits_{\tau\in\S_{m+n}^{(2)}}$ as 
$\sum\limits_{\tau\in\S\mn^{(2)}}+\sum\limits_{\tau\in\left(\S_{m+n}^{(2)}\setminus\S\mn^{(2)}\right)}$ 
and noticing that every 
$\tau\in\left(\S_{m+n}^{(2)}\setminus\S\mn^{(2)}\right)$ decomposes as the 
product of its restrictions $\tau_1$ and $\tau_2$ to the sets 
$\{1,2,\ldots,m\}$ and $\{m+1,m+2,\ldots,m+n\}$, it follows from part (1) 
of Lemma \ref{lemma2} and from part (1) of Lemma \ref{lemma1} that \begin{multline}\label{var}
\E\left(\tr\N\alpha\N\cdot\tr\N\beta\N\right)-\E\left(\tr\N\alpha\N\right)\cdot\E\left(\tr\N\beta\N\right)=
\\=\sum_{\substack{\tau\in\S\mn^{(2)}\text{ such that}\\r_a=r_{\tau(a)}\text{ and }\varepsilon_a\neq\varepsilon_{\tau(a)}\;\forall a\in[m+n]}}
     \frac{N^{\#(\tau)}}{N^{m+n+2}\cdot(N!)^s}\sum_{\sigma_1,\ldots,\sigma_s\in\S_N}
     \left(\prod_{C\in\tau^{-1}\gamma\mn}\Fix w_C(\sigma_1,\ldots,\sigma_s)\right)
\\+\sum_{\substack{\tau_1\in\S_m^{(2)}\text{ and }\tau_2\in\S_n^{(2)}\text{ such that}\\r_a=r_{\tau(a)}\text{ and }\varepsilon_a\neq\varepsilon_{\tau(a)}\;\forall a\in[m]
                                                               \\\text{and such that}\\r_{m+a}=r_{m+\tau(a)}\text{ and }\varepsilon_{m+a}\neq\varepsilon_{m+\tau(a)}\;\forall a\in[n]}}
                                \frac{N^{\#(\tau_1)+\#(\tau_2)}}{N^{m+n+2}}\Bigg[f_N(\tau_1,\tau_2)-g_N(\tau_1,\tau_2)\Bigg],
\end{multline}
where
\begin{equation*}
f_N(\tau_1,\tau_2):=\frac{1}{(N!)^s}\sum_{\sigma_1,\ldots,\sigma_s\in\S_N}
                                    \left(\prod_{\substack{C\in\tau_1^{-1}\gamma_m\\D\in\tau_2^{-1}\gamma_n}}
                                           \Fix w_C(\sigma_1,\ldots,\sigma_s)\cdot\Fix w_{m+D}(\sigma_1,\ldots,\sigma_s)\right)
\end{equation*}
and
\begin{multline*}
g_N(\tau_1,\tau_2):=\left[\frac{1}{(N!)^s}\sum_{\sigma_1,\ldots,\sigma_s\in\S_N}\left(\prod_{C\in\tau_1^{-1}\gamma_m}\Fix w_C(\sigma_1,\ldots,\sigma_s)\right)\right]\cdot
\\\cdot\left[\frac{1}{(N!)^s}\sum_{\sigma_1,\ldots,\sigma_s\in\S_N}\left(\prod_{D\in\tau_2^{-1}\gamma_n}\Fix w_{m+D}(\sigma_1,\ldots,\sigma_s)\right)\right].
\end{multline*}
\big(In the above, if  $\tau_2\in\S_n^{(2)}$ and if $D:=(a_1,a_2,\ldots,a_p)\in\tau_2^{-1}\gamma_n$, we denote
$$m+D:=(m+a_1,m+a_2,\ldots,m+a_p).\big)$$

To prove \eqref{alphabeta}, it now suffices to show that each term of the two sums appearing in \eqref{var} is $\mO(1/N^2)$.

If $\tau\in\S\mn^{(2)}$, then
\begin{align*}
&\frac{N^{\#(\tau)}}{N^{m+n+2}\cdot(N!)^s}\sum_{\sigma_1,\ldots,\sigma_s\in\S_N}\left(\prod_{C\in\tau^{-1}\gamma\mn}\Fix w_C(\sigma_1,\ldots,\sigma_s)\right)\leq
\\&\leq\frac{N^{\#(\tau)}}{N^{m+n+2}\cdot(N!)^s}\sum_{\sigma_1,\ldots,\sigma_s\in\S_N}N^{\#(\tau^{-1}\gamma\mn)}
=\frac{N^{\#(\tau)+\#(\tau^{-1}\gamma\mn)-(m+n)}}{N^2}
=\mO(1/N^2)\quad\text{by \eqref{nca}.}
\end{align*}
If $\tau_1\in\S_m^{(2)}$ is such that $w_C=e$ for every $C\in\tau^{-1}\gamma_m$, then it is easily seen that $f_N(\tau_1,\tau_2)=g_N(\tau_1,\tau_2)$.
Similarly, if $\tau_2\in\S_n^{(2)}$ is such that $w_{m+D}=e$ for every $D\in\tau_2^{-1}\gamma_n$, then $f_N(\tau_1,\tau_2)=g_N(\tau_1,\tau_2)$.

Finally, suppose that $\tau_1\in\S_m^{(2)}$ and $\tau_2\in\S_n^{(2)}$ are 
such that there exist $C_0\in\tau_1^{-1}\gamma_m$ and 
$D_0\in\tau_2^{-1}\gamma_n$ with $w_{C_0}\neq e\neq w_{m+D_0}$. Then
\begin{align*}
&\frac{N^{\#(\tau_1)+\#(\tau_2)}}{N^{m+n+2}}\cdot f_N(\tau_1,\tau_2)
\leq\frac{N^{\#(\tau_1)+\#(\tau_2)}}{N^{m+n+2}}\cdot\frac{1}{(N!)^s}\sum_{\sigma_1,\ldots,\sigma_s\in\S_N}N^{\#(\tau_1^{-1}\gamma_m)-1+\#(\tau_2^{-1}\gamma_n)-1}\cdot
                                                                                          \\&\hspace{2cm}\cdot\Fix w_{C_0}(\sigma_1\ldots,\sigma_s)\cdot\Fix w_{m+D_0}(\sigma_1\ldots,\sigma_s)=
\\&=\frac{N^{\#(\tau_1)+\#(\tau_1^{-1}\gamma_m)-(m+1)+\#(\tau_2)+\#(\tau_2^{-1}\gamma_n)-(n+1)}}{N^2}\cdot
\\&\hspace{2cm}\cdot\frac{1}{(N!)^s}\sum_{\sigma_1,\ldots,\sigma_s\in\S_N}\Fix w_{C_0}(\sigma_1\ldots,\sigma_s)\cdot\Fix w_{m+D_0}(\sigma_1,\ldots,\sigma_s)
=\mO(1/N^2)
\end{align*}
by \eqref{nc} and part (2) of Lemma \ref{Nica}, and
\begin{align*}
&\frac{N^{\#(\tau_1)+\#(\tau_2)}}{N^{m+n+2}}\cdot g_N(\tau_1,\tau_2)\leq
\\&\leq\frac{N^{\#(\tau_1)+\#(\tau_2)}}{N^{m+n+2}}
     \left(\frac{1}{(N!)^s}\sum_{\sigma_1,\ldots,\sigma_s\in\S_N}N^{\#(\tau_1^{-1}\gamma_m)-1}\cdot\Fix w_{C_0}(\sigma_1,\ldots,\sigma_s)\right)\cdot
     \\&\hspace{2cm}\cdot\left(\frac{1}{(N!)^s}\sum_{\sigma_1,\ldots,\sigma_s\in\S_N}N^{\#(\tau^{-1}\gamma_n)-1}\cdot\Fix w_{m+D_0}(\sigma_1,\ldots,\sigma_s)\right)=
\\&=\frac{N^{\#(\tau_1)+\#(\tau_1^{-1}\gamma_m)-(m+1)+\#(\tau_2)+\#(\tau_2^{-1}\gamma_n)-(n+1)}}{N^2}
       \left(\frac{1}{(N!)^s}\sum_{\sigma_1,\ldots,\sigma_s\in\S_N}\Fix w_{C_0}(\sigma_1,\sigma_2,\ldots,\sigma_s)\right)\cdot
       \\&\hspace{2cm}\cdot\left(\frac{1}{(N!)^s}\sum_{\sigma_1,\ldots,\sigma_s\in\S_N}\Fix w_{m+D_0}(\sigma_1,\sigma_2,\ldots,\sigma_s)\right)
=\mO(1/N^2)
\end{align*}
by \eqref{nc} and part (1) of Lemma \ref{Nica}.
\end{proof}
%%%
\begin{rem}\label{MN2}
The same method used to prove Theorem \ref{thm2} can be used to prove that, with the notations of Remark \ref{MN1}, we have that
\begin{equation}\label{HTUas}
\left(T_r\MN,U_r\MN,H_r\MN\right)_{r=1}^s\convergesas(T_r,U_r,H_r)_{r=1}^s
\end{equation}
as $M,N\to\infty$ with $M/N\to c$.

The main step in the proof is the following statement (analogous to Lemma \ref{lemma2}), whose proof we leave to the reader:

\emph{Let $k,l\geq1$, let $r_1,r_2,\ldots,r_{2k+2l}\in[s]$, and let $w_1,w_2,\ldots,w_{2k+2l}\in\F_s$. If 
$$\alpha=H_{r_1}^*T_{w_1}H_{r_2}U_{w_2}\cdots H_{r_{2k-1}}^*T_{w_{2k-1}}H_{r_{2k}}U_{w_{2k}}$$
and
$$\beta=H_{r_{2k+1}}^*T_{w_{2k+1}}H_{r_{2k+2}}U_{w_{2k+2}}\cdots H_{r_{2k+2l-1}}^*T_{w_{2k+2l-1}}H_{r_{2k+2l}}U_{w_{2k+2l}},$$
then for every $M,N\geq1$ we have that}
\begin{multline}\label{lemma2MN}
\E\left(\tr^{(M+N)}\alpha\MN\cdot\tr^{(M+N)}\beta\MN\right)
=\hspace{-0.5cm}\sum_{\substack{\tau\in\S_{2k+2l}^{(2)}\text{ such that}\\r_a=r_{\tau(a)}\;\forall a\in[2k+2l]}}
   \hspace{-0.5cm}\frac{(M+N)^{\#(\tau)}}{(M+N)^{2k+2l+2}\cdot(M!)^s(N!)^s}\cdot
\\\cdot\sum_{\substack{\sigma_1,\ldots,\sigma_s\in\S_M\\\rho_1,\ldots,\rho_s\in\S_N}}\left(\prod_{C\in(\tau^{-1}\gamma_{2k,2l})_{\mathrm{odd}}}\Fix w_C(\sigma_1,\ldots,\sigma_s)\right)
                                                                                                                                                         \left(\prod_{C\in(\tau^{-1}\gamma_{2k,2l})_{\mathrm{even}}}\Fix w_C(\rho_1,\ldots,\rho_s)\right).
\end{multline}
>From \eqref{lemma2MN} and \eqref{lemma1MN}, it follows by the same method used in the proof of Theorem \ref{thm2} that, if $\alpha$ and $\beta$ are as above and if $(M_q)_{q=1}^\infty$ and $(N_q)_{q=1}^\infty$ are two increasing sequences with $\lim\limits_{q\to\infty}\frac{M_q}{N_q}=c$, then
\begin{multline*}
\hspace{-2pt}\E\left(\tr^{(M_q+N_q)}\alpha\MNq\cdot\tr^{(M_q+N_q)}\beta\MNq\right)-\E\left(\tr^{(M_q+N_q)}\alpha\MNq\right)\cdot\E\left(\tr^{(M_q+N_q)}\beta\MNq\right)=
\\=\mO\left(1/N_q^2\right),
\end{multline*}
which implies that
\begin{equation*}
\sum_{q=1}^\infty\left[\E\left(\left|\tr^{(M_q+N_q)}\alpha\MNq\right|^2\right)-\left|\E\left(\tr^{(M_q+N_q)}\alpha\MNq\right)\right|^2\right]<\infty.
\end{equation*}
Now \eqref{HTUas} follows by an application of Lemma \ref{Steen}.
\end{rem}
%%%

\emph{Department of Pure Mathematics\\
University of Waterloo\\
200 University Avenue West\\
Waterloo, Ontario, Canada\\
N2L 3G1\\
e-mail: mgneagu@math.uwaterloo.ca}
%%%
\end{document}